\documentclass[11pt]{article} 
\usepackage{amsmath,amsthm} 
\usepackage{amssymb,latexsym}
\usepackage[T1]{fontenc}
\usepackage{authblk}
\usepackage{lipsum}
\usepackage{hyperref}
\usepackage{lipsum} % just for generating dummy text
\usepackage{lineno}
\usepackage{upgreek}
\usepackage{marvosym}
\modulolinenumbers[1]
\begin{document}
	\newtheorem{theorem}{Theorem}[section]
	\newtheorem{question}{Question}[section]
	\newtheorem{thm}[theorem]{Theorem}
	\newtheorem{lem}[theorem]{Lemma}
	\newtheorem{eg}[theorem]{Example}
	\newtheorem{prop}[theorem]{Proposition}
	\newtheorem{cor}[theorem]{Corollary}
	\newtheorem{rem}[theorem]{Remark}
	\newtheorem{deff}[theorem]{Definition}
	\numberwithin{equation}{section}
	\title{On $S$-$J$-Noetherian Rings}
	
	\author[1]{Tushar Singh}
	\author[2]{Ajim Uddin Ansari}
	\author[3]{Shiv Datt Kumar }

	\affil[1, 3]{\small Department of Mathematics, Motilal Nehru National Institute of Technology Allahabad, Prayagraj 211004, India \vskip0.01in Emails: sjstusharsingh0019@gmail.com, tushar.2021rma11@mnnit.ac.in, sdt@mnnit.ac.in}
	\vskip0.05in
	\affil[2]{\small Department of Mathematics, CMP Degree College, University of Allahabad, Prayagraj-211002, India \vskip0.01in Email: ajimmatau@gmail.com}
	\maketitle
	\hrule
	
	\begin{abstract}
		Let $R$ be a commutative ring with identity, $S\subseteq R$ be a multiplicative set and $J$ be an ideal of $R$. In this paper, we introduce the concept of $S$-$J$-Noetherian rings, which generalizes both $J$-Noetherian rings and $S$-Noetherian rings. We study several properties and charaterizations of this new class of rings. For instance, we prove Cohen's-type theorem for $S$-$J$-Noetherian rings. Among other results, we establish the existence of $S$-primary decomposition in $S$-$J$-Noetherian rings as a generalization of classical Lasker-Noether theorem.
	\end{abstract}
\hrule
\section{Introduction}
Throughout the paper, let $R$ be a commutative ring with identity, $S \subseteq R$ be a multiplicative set, and $J$ be a fixed ideal of $R$. For an ideal $I$ of $R$, we denote $\overline{S}=\{s+I\mid s\in S\}$ which is a multiplicative closed subset of $R/I$. The Noetherian property of rings plays a crucial role in areas such as commutative algebra and algebraic geometry. Given the significance of Noetherian rings, numerous authors attempted to generalize the concept of Noetherian rings (see \cite{ka24}, \cite{ad02}, \cite{ad24}, \cite{hk21}, \cite{jw14}, \cite{er22}, and  \cite{ts23}). As one of its crucial generilizations, Anderson and  Dumitrescu \cite{ad02} introduced the concept of  $S$-Noetherian rings. An ideal $I$ of $R$ is $S$-finite if there exists an element $s \in S$ and a finitely generated ideal $F$ of $R$ such that $sI \subseteq F \subseteq I$. A ring $R$ is called $S$-Noetherian if every ideal of $R$ is $S$-finite. Recently, Alhazmy et al. \cite{ka24} introduced the concept of $J$-Noetherian rings as a generalization of Noetherian ring. An ideal $I$ of $R$ is called a $J$-ideal if $I \nsubseteq J$ and $R$ is said to be $J$-Noetherian if every $J$-ideal is finitely generated. A particular interesting case occurs when $J=Nil(R)$, the ideal consisting of all nilpotent elements of $R$. In this situation, a $J$-Noetherian ring is referred to as a Nonnil-Noetherian ring, which was first introduced and studied by Badawi in \cite{ab03}. Furthermore, when $J=J(R)$, the Jacobson radical of $R$, the $J$-Noetherian ring is termed a non-$J$-Noetherian ring. This class of rings was first introduced by Dabbabi \cite{ad24} et al. in 2024, where they characterized various properties of non-$J$-Noetherian rings.

The primary objective of this paper is to introduce and study the notion of $S$-$J$-Noetherian rings. We present an example of an $S$-$J$-Noetherian ring which is not an $S$-Noetherian ring (see Example \ref{j9}). We generalize various properties and characterizations of both $J$-Noetherian and $S$-Noetherian rings to this new class of rings. For instance, we establish Cohen-type theorem for $S$-$J$-Noetherian rings and prove that the polynomial ring $R[X]$ is $S$-$J$-Noetherian if and only if it is $S$-Noetherian. Also, we show that the quotient of an $S$-$J$-Noetherian ring is an $\overline{S}$-Noetherian ring (see Proposition \ref{j6}). Moreover, we provide necessary and sufficient conditions for an $S$-$J$-Noetherian ring to belong to the class of $S$-Noetherian rings (see Theorem \ref{j3} and \ref{j10}). In \cite[Theorem 2.10]{ts24}, among the other result, Singh et al. generalized the classical Lasker-Noether theorem for $S$-Noetherian modules. We end the paper by extending the classical Lasker-Noether theorem for the class of $S$-$J$-Noetherian rings (see Theorem~\ref{on}).

	\section{Main Results}
	We begin by introducing the concept of $S$-$J$-Noetherian rings. 
	\begin{deff}
			Let $R$ be a ring, $S\subseteq R$ be a multiplicative set, and $J$ an ideal of $R$. An ideal $I$ of $R$ is said to be a $J$-ideal if $I \nsubseteq J$. We say that $R$ is an $S$-$J$-Noetherian ring if each $J$-ideal of $R$ is $S$-finite.
	\end{deff}
It is evident that every $J$-Noetherian ring is an $S$-$J$-Noetherian ring when $S=\{1\}$. However, the following example illustrates that the converse is not true in general.

\begin{eg}\label{j8}
	\noindent
Consider the ring $R = \mathcal{F}[X_1, X_2, \dots]$, where $\mathcal{F}$ is a field, and let $J = (0)$. Define the ideal $I = (X_1, \dots, X_{n}, \ldots)$. Clearly, $I$ is a $J$-ideal but is not finitely generated. Hence $R$ is not a $J$-Noetherian ring. Now, let $S = R \setminus \{0\}$ be the multiplicative closed subset of $R$. Let $K$ be a nonzero proper ideal of $R$. Evidently, $K$ is $J$-ideal and $K\cap S\neq\emptyset$. Therefore, by \cite[Proposition 2(a)]{ad02}, $K$ is $S$-finite. Hence $R$ is $S$-$J$-Noetherian.  
\end{eg}

Cohen's theorem is the classical result which states that a ring is Noetherian if all its prime ideals are finitely generated. Now, we extend this result for $S$-$J$-Noetherian rings.

\begin{theorem}\label{j4}
	A ring $R$ is $S$-$J$-Noetherian if and only if its prime $J$-ideals (disjoint from $S$) are $S$-finite.
\end{theorem}
\begin{proof}
	If $R$ is $S$-$J$-Noetherian, then it is obvious that all prime $J$-ideals of $R$ are $S$-finite. Now, suppose that all prime $J$-ideals (disjoint from $S$) of $R$ are $S$-finite and assume that $R$ is not $S$-$J$-Noetherian. Therefore the set $\mathcal{F}$ of all $J$-ideals that are non-$S$-finite is a non-empty set which is ordered by the inclusion. By Zorn's lemma, choose $P$ maximal in $\mathcal{F}$. This implies $P$ is not a $S$-finite and so $P\cap S=\emptyset$. We show that $P$ is a prime ideal of $R$. This makes $P$ a $J$-prime ideal (disjoint from $S$) that is $S$-finite, which is a contradiction to the fact that $P \in \mathcal{F}$. Suppose there exists $a, b \in R\setminus P$ such that $a b \in P$. If $P+aR\subseteq J$, $P\subseteq J$, a contradiction, as $P$ is $J$-ideal. Therefore $P+aR$ is $J$-ideal.
	Since $P\subsetneq P+aR$, it follows that $P+aR$ is $S$-finite since $P$ is a maximal element of $\mathcal{F}$. Then, there exist $s\in S$, $\alpha_{1},\ldots, \alpha_{n}\in P$ and $x_{1},\ldots x_{n}\in R$ such that $s(P+aR)\subseteq \left(\alpha_{1}+a x_{1}, \ldots,  \alpha_{n}+a x_{n}\right)\subseteq P+aR$. Consider the ideal $Q=(P:a) =\{x\in R\mid ax\in P\}$. Evidently, $Q$ is $J$-ideal and $P\subsetneq Q$ as $b\in Q\setminus P$. By the maximalty of $P$, $Q$ is an $S$-finite ideal. Then there exist $t\in S$ and $\beta_ {1}, \ldots, \beta_{k} \in Q$ such that $tQ\subseteq \left(\beta_{1}, \ldots, \beta_{k}\right)\subseteq Q$. Let $x\in P$. Then $sx\in s(P+aR)\subseteq\left(\alpha_{1}+a x_{1}, \ldots, \alpha_{n}+a x_{n}\right)$, and so there exist $u_{1}, \ldots, u_{n} \in R$ such that $sx=u_{1}\left(\alpha_{1}+a x_{1}\right)+\cdots+u_{n}\left(\alpha_{n}+a x_{n}\right)=u_{1} \alpha_{1}+\cdots+u_{n} a_{n}+a\left(u_{1} x_{1}+\cdots+u_{n} x_{n}\right)$. So $a\left(u_{1} x_{1}+\cdots+u_{n} x_{n}\right) =sx-\left(u_{1} \alpha_{1}+\cdots+u_{n}\alpha_{n}\right) \in P$. Then $u_{1} x_{1}+\cdots+u_{n} x_{n}\in (P:a)=Q$. Therefore we can find $w_{1},\ldots, w_{k}\in R$ such that $t\left(u_{1} x_{1}+\ldots+u_{n} x_{n}\right)=w_{1}\beta_{1}+\cdots+w_{k}\beta_{k}$, which states that $stx=t\left(u_{1}\alpha_{1}+\cdots+u_{n} a_{n}\right)+at\left(u_{1} x_{1}+\cdots+u_{n} x_{n}\right)=t\left(u_{1}\alpha_{1}+\cdots+u_{n} a_{n}\right)+a \left( w_{1}\beta_{1}+\cdots+w_{k}\beta_{k}\right)$. Hence we obtain $uP\subseteq \left(\alpha_{1}, \ldots, \alpha_{n}, a \beta_{1}, \ldots, a \beta_{k}\right)\subseteq P$, where $u=st\in S$, which means that $P$ is $S$-finite. This contradicts to the choice of $P$. Thus $R$ is $S$-$J$-Noetherian.
\end{proof}

Every $S$-Noetherian ring is clearly an $S$-$J$-Noetherian ring. However, an $S$-$J$-Noetherian ring need not be an $S$-Noetherian ring. For this, consider the following example.
\begin{eg}\label{j9}
Consider a ring $R_1=\mathcal{F}[X_1, \ldots, X_{n}, \ldots]$, where $\mathcal{F}$ is a field, and $I=(X^{2}_{i}; i\in\mathbb{N})$ be an ideal of $R_1$. Let $R=R_1/I$. Consider the prime ideal $P = (X_i; i\in\mathbb{N})$ of $R_1$. Note that any prime ideal of the ring $R$ contains $P/I$. Then the unique minimal prime ideal of $R$ is $P /I$. Take $J=P/I$ and $S=R\setminus (P/I)$ is a multiplicative subset of $R$. Any $J$-prime ideal $P'$ of $R$ contains properly $P/I$, and then $P'\cap S\neq\emptyset$. By \cite[Proposition 2(a)]{ad02}, $P'$ is $S$-finite. By Theorem \ref{j4}, $R$ is an $S$-$J$-Noetherian ring. Next, our aim is to show that $R$ is not a $S$-Noetherian ring. Suppose the ideal $P /I$ is $S$-finite. There exist $\bar{s}\in S$ and $i_1, \dots, i_n \in \mathbb{N}$ such that $\bar{s}(P /I) \subseteq (\overline{X_{i_1}}, \dots, \overline{X_{i_n}}) \subseteq (P /I)$. The polynomial $\bar{s}$ of $R$ uses a finite number of variables $X_{j_1}, \dots, X_{j_m}$ and its constant term $d \neq 0$. Let $k \in \mathbb{N} \setminus \{i_1, \dots, i_n, j_1, \dots, j_m\}$. Then, $\bar{s}\overline{X_k} = \bar{f}_1 \overline{X_{i_1}} + \dots + \bar{f}_n \overline{X_{i_n}}$, where $\bar{f}_1, \dots, \bar{f}_n \in R$. Thus $sX_k- f_1X_{i_1} - \dots -f_n X_{i_n}\in I$. This implies that $X_{i_1} = \dots = X_{i_n} = X_{j_1} = \dots = X_{j_m} = 0$, we obtain $d X_k \in \left( X_i^2 \mid i \in \mathbb{N} \setminus \{i_1, \dots, i_n, j_1, \dots, j_m\} \right)$. This is a contradiction.
\end{eg}
Examples \ref{j8} and \ref{j9} demonstrate that the concept of  \( S \)-\( J \)-Noetherian rings is a proper generalization of both the \( J \)-Noetherian rings and  \( S \)-Noetherian rings. 

Recall \cite{zb17}, let $E$ be a family of ideals of a ring $R$. An element $I\in E$ is said to be an \textit{$S$-maximal element} of $E$ if there exists an $s\in S$ such that for each $J\in E$, if $I\subseteq J$, then $sJ\subseteq I$. Also, a chain of ideals $(I_{i})_{i\in\wedge}$  of $R$ is called \textit{$S$-stationary} if there exist $k\in \wedge$ and $s\in S$ such that  $sI_{i}\subseteq I_{k}$ for all $i\in\wedge$, where $\wedge$ is an arbitrary indexing set.  A family $\mathcal{F}$ of ideals of $R$ is said to be $S$-saturated if it satisfies the following property: for every ideal $I$ of $R$, if there exist $s\in S$ and $J\in\mathcal{F}$ such that $sI\subseteq J$, then $I\in\mathcal{F}$.

\begin{theorem}\label{j1}
	Let $J$ be a proper ideal of $R$. Then the following statements are equivalent.
	\begin{enumerate}
		\item $R$ is an $S$-$J$-Noetherian.
		
		\item  Every ascending chain of $J$-ideals of $R$ is $S$-stationary.
		
		\item Every nonempty $S$-saturated set of  $J$-ideals of $R$ has a maximal element.
		
		\item Every nonempty family of $J$-ideals has an $S$-maximal element with respect to inclusion.
	\end{enumerate}
\end{theorem}

\begin{proof}
	\leavevmode
	
		\item(1) $\Rightarrow$ (2). Let $(I_n)_{n\in\wedge}$ be an increasing sequence of $J$-ideals of $R$. Define the ideal $I = \bigcup\limits_{n\in\wedge}I_n$. If $I\subseteq J$, then $I_{n}\subseteq J$, which is not possible since each $I_n$ is a $J$-ideal. Thus $I$ is a  $J$-ideal of $R$. Also, $I$ is $S$-finite since $R$ is $S$-$J$-Noetherian. Consequently, there exist a finitely  generated ideal $F\subseteq R$ and $s\in S$ such that $sI\subseteq F\subseteq I$. Since $F$ is finitely generated, there is a $k\in\wedge$ satisfying $F\subseteq I_{k}$. Then we have $sI\subseteq F\subseteq I_{k}$, from which it follows that $sI_n \subseteq I_{k}$ for each $n\in\wedge$.\\
		
		\item (2) $\Rightarrow$ (3).  Let $\mathcal{D}$ be an $S$-saturated set of $J$-ideals of $R$. Given any chain \(\{I_n\}_{n \in\wedge} \subseteq \mathcal{D} \), we claim that \( I = \bigcup\limits_{n\in\wedge} I_n \) belongs to $\mathcal{D}$, which will establish that \( I \) as an upper bound for the chain. Indeed, by (2), there exist \( k \in \wedge \) and \( s \in S \) such that \( s I_n \subseteq I_k \) for every \( n \in \wedge \). Consequently, we obtain $sI = s \left(\bigcup\limits_{n\in\wedge} I_n \right) \subseteq I_k$. Since $\mathcal{D}$ is $S$-saturated, it follows that \( I \in \mathcal{D} \), as required. Applying Zorn's lemma, we conclude that $\mathcal{D}$ has a maximal element.\\
		
		\item (3) $\Rightarrow$ (4). Let $\mathcal{D}$ be a nonempty set of  $J$-ideals of $R$. Consider the family $\mathcal{D}^{S}$ of all $J$-ideals $L \subseteq R$ such that there exist some $s \in S$ and $L_0 \in \mathcal{D}$ with $sL \subseteq L_0$. Clearly, $\mathcal{D}\subseteq \mathcal{D}^{S}$, so $\mathcal{D}^{S}\neq\emptyset$. It is straightforward to see that $\mathcal{D}^{S}$ is $S$-saturated. Thus, by (3) $\mathcal{D}^{S}$ has a maximal element $K\in\mathcal{D}^{S}$. Fix $s \in S$ and $Q\in \mathcal{D}$ such that $sK\subseteq Q$. Now, we claim that $Q$ is an $S$-maximal element of $\mathcal{D}$; specifically, given $L\in \mathcal{D}$ with $Q\subseteq L$, we will show that $sL \subseteq Q$. Note that $K + L$ satisfies $s(K + L)= sK + sL \subseteq Q + L \subseteq L$, so that $K + L \in \mathcal{D}^{S}$. Also, if $(K+L)\subseteq J$, then $K\subseteq J$, which is not possible since $K$ is a $J$-ideal of $R$. Thus $K+L$ is a $J$-ideal of $R$. Therefore maximality of $K$ implies $K = K + L$, so that $L \subseteq K$. But then $sL \subseteq sK \subseteq Q$, as desired.\\
		
		\item (4) $\Rightarrow$ (1). Let $I$ be a $J$-ideal of $R$, which we will prove to be $S$-finite. Let $\mathcal{D}$ be the family of finitely generated $J$-ideal of $R$ such that $J\subseteq I$. Choose $x\in I\setminus J$. Then $L=(x)\subseteq I$, and $L\nsubseteq J$. This implies that $L\in\mathcal{D}$, and so $\mathcal{D}$ is nonempty. Then $\mathcal{D}$ has an $S$-maximal element $K\in\mathcal{D}$. Fixing $x \in I$, take a finitely generated ideal of the form $Q= K + xR$. Since $K\subseteq I$ and $x\in I$, so $Q\subseteq I$. Consequently, $Q\in\mathcal{D}$ such that $K\subseteq Q$. This implies that there exists $s\in S$ such that $sQ\subseteq K$; in particular, $sx\in K$. This verifies $sI\subseteq K\subseteq I$, so that $I$ is $S$-finite. It follows that $R$ is $S$-$J$-Noetherian.

\end{proof}

Let $f : R \to R'$ be a homomorphism and $S$ a multiplicative closed subset of $R$. Then it is easy to see that $f(S)$ is a multiplicative closed subset of $R'$ if $0 \notin f(S)$
and $1 \in f(S)$.
\begin{prop}\label{j12}
	Let $f:R\to R'$ be an epimorphism and $J$ be an ideal of $R'$. If $R$ is an $S$-$f^{-1}(J)$-Noetherian ring with $0 \notin f(S)$, then $R'$ is a $f(S)$-$J$-Noetherian ring. 
\end{prop}
\begin{proof}
	 Suppose  $\{I_{i}\} _{i\in\wedge}$ is any increasing chain of $J$-ideals of $R'$. Then $I_{i}\nsubseteq J$ for each $i\in\wedge$. Suppose contrary that,  for each $i$ there exist $\alpha_{i}\in I_{i}\setminus J$ such that $f^{-1}(\alpha_{i})\subseteq f^{-1}(J)$. Then $\alpha_{i}\in f(f^{-1}(\alpha_{i}))\subseteq f(f^{-1}(J))=J$, for $f$ is an epimorphism. This is a contradiction, as $\alpha_{i}\notin J$. Thus $f^{-1}(I_{i})\nsubseteq f^{-1}(J)$ for each $i\in\wedge$ and hence  $f^{-1}(I_{i})$ is  $f^{-1}(J)$  ideal of $R$. Then we have an increasing chain $\{f^{-1}(I_{i})\}_{i\in\wedge}$ of $f^{-1}(J)$-ideal of $R$. Since $R$ is an $S$-$f^{-1}(J)$-Noetherian, there exist $k\in\wedge$ and $s\in S$ such that $sf^{-1}(I_i)\subseteq f^{-1}(I_k)$ for all $i\in\wedge$. Applying $f$ to both sides, we obtain $f(s f^{-1}(I_i))= f(s) f(f^{-1}(I_i))\subseteq f(f^{-1}(I_k))$ for all $i\in\wedge$. Since $f$ is an epimorphism, it follows that $f(s) I_i \subseteq I_k$ for all $i\in\wedge$. Hence, by Theorem \ref{j1}, $R'$ is a $f(S)$-$J$-Noetherian ring. 
\end{proof}

\begin{theorem}\label{j3}
	Let $S$ be a multiplicative subset of a ring $R$. The following
	statements are equivalent:
	\begin{enumerate}
		\item $R$ is S-Noetherian.
		\item $R$ is $S$-$J$-Noetherian and $J$ is an $S$-finite ideal of $R$.
	\end{enumerate}
\end{theorem}
\begin{proof}
	(1) $\Rightarrow$ (2). This implication is obvious. (2) $\Rightarrow$ (1). Let $P$ be a prime ideal of $R$. If $P\subseteq J$, then $P$ is $S$-finite by the assumption. Suppose that $P$ contains properly in $J$. Then $P$ is a $J$-ideal of $R$ disjoint with $S$. Since $R$ is $S$-$J$-Noetherian, then $P$ is $S$-finite disjoint from $S$. So, by \cite[Corollary 5]{ad02}, $R$ is $S$-Noetherian.
\end{proof}
Let $R$ be a ring and $S$ be a multiplicative subset of $R$. Recall \cite{ad02}, let $S$ be an anti-Archimedean subset of $R$ if $\bigcap_{n \geq 1} s^n R \cap S \neq \emptyset$ for all  $s \in S$. 
\begin{cor}
	Let $S\subseteq R$ be an anti-Archimedean multiplicative set and $J$ is $S$-finite. If $R$ is $S$-$J$-Noetherian, then the polynomial ring $R[X_1,\ldots, X_n]$ is also $S$-$J$-Noetherian.
\end{cor}
\begin{proof}
	By Theorem \ref{j3}, $R$ is $S$-Noetherian ring. Then, by \cite[Proposition 9]{ad02}, $R[X_1,\ldots, X_n]$ is $S$-Noetherian. This implies $R[X_1,\ldots, X_n]$ is $S$-$J$-Noetherian.
\end{proof}
Recall \cite{ad09}, let $M$ be an $R$-module. The the idealization of $R$-module $M$, $R(+)M=\{(r, m)\mid r\in R, m\in M\}$ is a commutative ring with componentwise addition and  multiplication defined by $(\alpha_1, m_1)(\alpha_2, m_2)= (\alpha_1\alpha_2, \alpha_1m_2 + \alpha_2m_1)$ for all $\alpha_1, \alpha_2 \in R$ and $m_1, m_2 \in M$. It is straightforward to verify that $S(+) M = \{(s, m) \mid s \in S, m \in M\}$ forms a multiplicative set in $R(+)M$.
\noindent
The following example shows that the polynomial ring over an $S$-$J$-Noetherian ring need not be $S$-$J$-Noetherian.
\begin{eg}
	Let $V$ be an infinite dimensional vector space over a field $K$. Then $R= K(+)V$ is an $S$-$J$-Noetherian ring for every multiplicative subset $S$ of $R$. Moreover, if $0 \notin S$, then $R[X]$ is not an $S$-$J$-Noetherian ring. In particular, if $J=Nil(R)$, then the proof follows from \cite[Example 2.4]{nm24}.
\end{eg}

We next show that the polynomial ring $R[X]$ is $S$-$J$-Noetherian if and only if it is $S$-Noetherian.
\begin{cor}\label{j5}
	Let $R$ be a ring, $S\subseteq R$ be a multiplicative set and $J$ be an ideal of $R$. Then $R[X]$ is an $S$-$J[X]$-Noetherian ring if and only if $R[X]$ is an $S$-Noetherian ring.
\end{cor}
\begin{proof}
	Suppose $R[X]$ is an $S$-$J[X]$-Noetherian ring. Then we show that $R[X]$ is an $S$-Noetherian ring. To prove this, by Theorem \ref{j3}, it is sufficient to show that $J[X]$ is $S$-finite. Define the ideal $Q=J[X]+XR[X]$ of $R[X]$. Note that $Q$ is a $J[X]$-ideal since $Q \not\subseteq J[X]$. Therefore $Q$ is $S$-finite. So there exist $s\in S$ and $f_1,\ldots, f_n\in R[X]$ such that $s(J[X]+XR[X])\subseteq f_1R[X]+\cdots+f_nR[X]\subseteq J[X]+XR[X]$. As a result, we get $sJ\subseteq f_1(0)R+\cdots+f_n(0)R\subseteq J$. This implies that $sJ[X]\subseteq f_1(0)R[X]+\cdots+f_n(0)R[X]\subseteq J[X]$. Thus $J[X]$ is an $S$-finite ideal of $R[X]$. The converse is trivially true.
\end{proof}

\begin{prop}\label{j6}
	Let $R$ be an $S$-$J$-Noetherian ring. Then $R/J$ is an $\overline{S}$-Noetherian ring.
\end{prop} 
\begin{proof}
	A nonzero prime ideal (disjoint from $\overline{S}$) of $R / J$ is of the form $P / J$ with $P \in \operatorname{Spec}(R)$ and $J \subsetneq P$. Evidently, $P$ is a $J$-ideal with $P\cap S=\emptyset$ since $P/J$ is nonzero and $P/J\cap\overline{S}=\emptyset$. By the hypothesis, $P$ is $S$-finite. Then there exist $s\in S$ and $p_1,\ldots, p_{n}\in P$ such that $sP\subseteq (p_1,\ldots, p_n)\subseteq P$. Let $x\in P$. Then we can find $a_1,\ldots, a_n\in R$ such that $sx=a_1p_1+\cdots+a_np_n$. It follows that $(s+J)(x+J)=(a_1+J)(p_1+J)+\ldots+(a_n+J)(p_n+J)$, where $s+J\in\overline{S}$ and $a_{1}+J,\ldots, a_{n}+J\in R/J$. This implies that $(s+J)(P/J)\subseteq (p_1+J,\ldots, p_n+J)\subseteq P/J$, i.e., $P/J$ is $\overline{S}$-finite. By \cite[Corollary 5]{ad02}, $R /J$ is $\overline{S}$-Noetherian.
\end{proof}
\begin{cor}
	Let $S\subseteq R$ be an anti-Archimedean multiplicative set. If $R$ is $S$-$J$-Noetherian, then polynomial ring $(R/J)[X_1,\ldots, X_n]$ is $\overline{S}$-$J[X_1,\ldots, X_n]$-Noetherian.
\end{cor}
\begin{proof}
	By Proposition \ref{j6}, $R/J$ is $\overline{S}$-Noetherian. Then, by \cite[Proposition 9]{ad02}, $(R/J)[X_1,\ldots, X_n]$ is $\overline{S}$-Noetherian. This implies $(R/J)[X_1,\ldots, X_n]$ is $\overline{S}$-$J[X_1,\ldots, X_n]$-Noetherian.
\end{proof}
\begin{deff}\cite{ab03}
	An ideal $I$ of a ring $R$ is called divided if $I\subset x R$ for every $x \in R \backslash I$.
\end{deff}
\begin{theorem}
	Let $R$ be an $S$-$J$-Noetherian ring, and $I$ be a $J$-ideal of $R$ disjoint from $S$. If $J$ is divided ideal, then there exist $s\in S$ and $S$-prime ideals  $P_1, \ldots, P_n$ of $R$ such that $s\left(P_1 \cdots P_n\right)\subseteq I$.
\end{theorem}
\begin{proof}
	Since $I\nsubseteq J$ and $J$ is divided, then $J\subset (x)\subseteq I$ for some $x\in I\setminus J$. Thus $I/J$ is an ideal of the $\overline{S}$-Noetherian ring $R/J$. Since $I\cap S=\emptyset$, then $\left(I/J\right)\cap\overline{S}=\emptyset$. For this, if $\left(I/J\right)\cap\overline{S}\neq\emptyset$, then $s+J=i+J$ for some $s\in S$ and $i\in I$. Consequently,  $s-i\in J\subset I$, and so $s\in I$, a contradiction as $I\cap S=\emptyset$. Thus $I/J$ is disjoint from $\overline{S}$. It follows that there exist $\bar{s}\in\overline{S}$ and $\overline{S}$-prime ideals $Q_1, \dots, Q_n$ of $R/J$ containing $I/J$ such that $\bar{s}(Q_1 \dots Q_n )\subseteq I/J$, by \cite[Theorem 5]{ah20}. Clearly, $Q_{i}\cap\overline{S}=\emptyset$ for each $i=1,\ldots,n$ since each $Q_{i}$ is $\overline{S}$-prime.
	Then, by \cite[Proposition 3]{ah20}, for each $1 \leq i \leq n$, there exists an $S$-prime ideal $P_i$ of $R$ containing $J$ such that $Q_i = P_i / J$. Therefore $\bar{s}\left(\left (P_1 \cdots P_n \right)/J\right) \subseteq I/J$ since $P_1/J\cdots P_{n}/J=\left( P_1\cdots P_{n} \right)/J$. For every $a\in P_1\cdots P_{n}$, $(s+J)(a+J)=b+J$ for some $b\in I$. Conseuently, $sa-b\in J\subset I$, and so $sa\in I+Rb\subseteq I+J=I$. Thus $s\left(P_1 \cdots P_n\right)\subseteq I$.
\end{proof}
\begin{prop}
	Let $R\subseteq R'$ be an extension of rings such that $IR' \cap R= I$ for each ideal $I$ of $R$, and let $S \subseteq R$ be a multiplicative set. If $R'$ is an $S$-$J$-Noetherian ring, then  $R$ is $S$-$J$-Noetherian.
\end{prop}
\begin{proof}
	Let $I$ be a $J$-ideal of $R$ and $I\subseteq IR'$. If $IR'\subseteq J$, then $I\subseteq J$, which is not possible since $I\nsubseteq J$. Thus $IR'$ is a $J$-ideal of $R'$.	Since the ring $R'$ is $S$-$J$-Noetherian, there exist $s \in S$ and  
	$i_1, \dots, i_n \in I$ such that $sIR' \subseteq (i_1, \dots, i_n)R'\subseteq IR'$. By hypothesis, $sI = sIR' \cap R\subseteq (i_1, \dots, i_n)R'\cap R \subseteq IR'\cap R=I$. Then $I$ is an $S$-finite ideal of $R$, as desired.
\end{proof}

\begin{prop}
	Let $R$ be an $S$-$J$-Noetherian ring and $I$ be a $J$-ideal of $R$ disjoint from $S$. Then there exist $t\in S$ and $m\in\mathbb{N}$ such that $t(rad(I))^{m}\subseteq I$.
\end{prop}
\begin{proof}
	Let $I$ be a $J$-ideal of $R$. Then $rad(I)$ is also a $J$-ideal of $R$, and hence $rad(I)$ is $S$-finite. Consequently, there exist  $s\in S$ and $x_1, \dots, x_n\in rad(I)$ such that $s(rad(I))\subseteq K \subseteq rad(I)$, where $K=(x_1, \dots, x_n)$. Suppose $m_{i}\in\mathbb{N}$ be such that $x_i^{m_i}\in I$ for any $1 \leq i \leq n$. Then  choose sufficiently large $m\in\mathbb{N}$ such that $K^{m}\subseteq I$. Therefore $t(rad(I))^{m}\subseteq I$, where $t=s^{m}\in S$.
\end{proof}
\begin{lem}\label{j7}
	Let $R$ be an $S$-$J$-Noetherian and $I$ be an $J$-ideal of $R$. Then $R/I$ is an $\overline{S}$-Noetherian ring.
\end{lem}
\begin{proof}
	Let $\{I_i/I\}_{i\in\wedge}$ be an ascending chain of non-zero ideals of $R / I$. As a result, $\{I_{i}\}_{i\in\wedge}$ is an ascending chain of $J$-ideal of $R$ and hence, by Theorem \ref{j1}, there exist $s\in S$ and $k \in\wedge$ such that $sI_{i}\subseteq I_{k}$ for every $i\in\wedge$. Therefore $(s+I)(I_{i} / I)\subseteq I_{k} / I$ for every $i\in\wedge$ and hence $\left(I_{i} / I\right)_{n \in\wedge}$ is $\overline{S}$-stationary. By \cite[Theorem 2.3]{zb17},  $R / I$ is $\overline{S}$-Noetherian.
\end{proof}
Recall that a ring $R$ is said to be decomposable if $R$ admits a non-trivial idempotent. Let $Idem(R)$ denote the set of idempotent elements of $R$.
\begin{theorem}\label{j10}
	Let $R$ be a decomposable ring and $J$ be an ideal of $R$ with $eJ \neq (e)$ for each $e \in \text{Idem}(R) \setminus \{0,1\}$. Then $R$ is $S$-$J$-Noetherian if and only if $R$ is  $S$-Noetherian.
\end{theorem}
\begin{proof}
	It is sufficient to prove that if $R$ is $S$-$J$-Noetherian, then  $R$ is $S$-Noetherian. To prove this, first we prove that $R/(e)$ is $\overline{S}$-Noetherian for each $e\in Idem(R)\setminus\{0,1\}$. Consider $e\in Idem(R)\setminus\{0,1\}$. Let $L$ be an ideal of $R$ which contains $(e)$. Then $e\notin J$ since $eJ \neq (e)$, and so $L\nsubseteq J$. Thus $L$ is a $J$-ideal, and so by Lemma \ref{j7}, $R/L$ is $\overline{S}$-Noetherian. This implies that  $R/(e)$ is $\overline{S}$-Noetherian since $(e)\subseteq L$. Now, let $K$ be an ideal of $R$ such that $K\subseteq (e)$ for each $e \in Idem(R) \setminus \{0,1\}$. We claim that $K$ is $S$-finite. Clearly, $eK=K$. If $K=(0)$, then $K$ is $S$-finite. So we may assume that $K\neq 0$. If $K\subseteq (1-e)$, then $eK\subseteq (e-e^{2})=(0)$, i.e., $eK=K=0$, a contradiction as $K\neq 0$. Therefore $K\nsubseteq (1-e)$. Since $1-e\in Idem(R)\setminus \{0,1\}$, $R/(1-e)$ is a $\overline{S}$-Noetherian ring. Set $I=(1-e)$ for simplicity. Then $L=\left(K+I\right)/I$ is an $\overline{S}$-finite ideal of $R/I$. Then there exist $\alpha_1+I,\ldots, \alpha_n+I\in R/I$, where $\alpha_1,\ldots, \alpha_n\in K$ and $s'=s+I\in\overline{S}$ such that $s'L\subseteq \left(\alpha_1+I,\ldots, \alpha_n+I\right)\subseteq L$. Let $\beta\in K+I$. Then $\beta+I\in L$, and so $s\beta+I\in s'L\subseteq \left(\alpha_1+I,\ldots, \alpha_n+I\right)$. This implies that $s\beta+I=(u_{1}+I)(\alpha_1+I)+\cdots+(u_{n}+I)(\alpha_n+I)$ for some $u_{1}+I, \ldots, u_{n}+I\in R/I$. Consequently, $s\beta-(u_1\alpha_{1}+\cdots+u_n\alpha_{n})\in I$, $s\beta-(u_1\alpha_{1}+\cdots+u_n\alpha_{n})\in F$, where $F =(\alpha_{1},\ldots, \alpha_{n}, 1-e)$ since $I\subseteq F$. Thus $s\beta\in F$, and hence $s(K+(1-e))\subseteq F\subseteq K+(1-e)$. Therefore $K+(1-e)$ is $S$-finite. Consequently, $K=Ke=(K+(1-e))e$ is an $S$-finite ideal of $R$, as claimed. Now, let $T$ be an ideal of $R$. Since $eT\subseteq (e)$ and $(1 - e) T\subseteq K+(1-e)T\subseteq K+(1-e)$ for each $e\in Idem(R)\setminus\{0, 1\}$, $eT$ and $(1-e)T$ are $S$-finite. It follows that $T = eT+ (1 - e)T$ is $S$-finite, and hence $R$ is $S$-Noetherian ring.
\end{proof}

\begin{deff}\label{S-iir.}\cite{ts24}
	An ideal $Q$ (disjoint from $S$) of the ring $R$ is called  $S$-irreducible if $s(I\cap K)\subseteq Q \subseteq I\cap K$ for some $s\in S$ and some ideals $I$, $K$ of $R$, then there exists $s'\in S$ such that either $ss'I\subseteq Q$ or $ss'K\subseteq Q$.
\end{deff}

\noindent
It is clear from the definition that every irreducible ideal is an $S$-irreducible ideal. However, the following example shows that an $S$-irreducibile ideal need not be irreducible. 
	\begin{eg}\label{fm}
	\noindent
	Let $R=\mathbb{Z}$, $S=\mathbb{Z}\setminus 3\mathbb{Z}$ and $I=6\mathbb{Z}$. Since $I=2\mathbb{Z}\cap 3\mathbb{Z}$, therefore $I$ is not an irreducible ideal of $R$. Now, take $s=2\in S$. Then $2(3\mathbb{Z})=6\mathbb{Z}\subseteq I$. Thus $I$ is an $S$-irreducible ideal of $R$.
\end{eg}
\noindent
Recall \cite[Definition 2.1]{me22}, a proper ideal $Q$ of a ring $R$ disjoint from $S$ is said to be $S$-primary if there exists an $s \in S$ such that for all $a, b \in R$, if $ab \in Q$, then either $sa \in Q$ or $sb \in rad(Q)$. Following from \cite{ts24}, let
$I$ be an ideal of $R$ such that $I\cap S=\emptyset$. Then $I$ admits $S$-primary decomposition if $I$ can be written as a finite intersection of $S$-primary ideals of $R$.\\

\noindent
Now, we extend $S$-primary decomposition theorem for $S$-$J$-Noetherian rings. We start with the following lemma.
	\begin{lem}\label{sp}
	\noindent
	Let $R$ be an $S$-$J$-Noetherian ring. Then every $S$-irreducible $J$-ideal of $R$ is $S$-primary.
\end{lem}
\begin{proof}
	Suppose $Q$ is an $S$-irreducible $J$-ideal of $R$. Let $a,b\in R$ be such that $ab\in Q$ and $sb\notin Q$ for all $s\in S$. Our aim is to show that there exists $t\in S$ such that $ta\in rad(Q)$.
	Consider $A_n=\{x\in R\hspace{0.2cm}|\hspace{0.1cm} a^{n}x\in Q \}$ for $n\in\mathbb{N}$. Since $Q$ is a $J$ ideal, there exists $\alpha\in Q\setminus J$. Then $a^{n}\alpha\in Q$ for each $n\in\mathbb{N}$. This implies that $\alpha\in A_{n}$ but $\alpha\notin J$ for each $n\in\mathbb{N}$. Consequently, each $A_{n}$ is a $J$-ideal of $R$ and $A_1\subseteq A_2\subseteq A_3\subseteq \cdots$ is an increasing chain of ideals of $R$. Since $R$ is a $S$-$J$-Noetherian, by Theorem \ref{j1}, this chain is $S$-stationary, i.e., there exist $k\in \mathbb{N}$ and $s\in S$ such that $sA_n\subseteq A_k$ for all $n\geq k$. Consider the two ideals $I=(a^{k}) +\hspace{0.1cm} Q$ and $K=(b) +\hspace{0.1cm} Q$ of $R$. Then  $Q\subseteq I\cap K$. For the reverse containment, let $y\in I\cap K$. Write  $y=a^{k}z+q$ for some $z\in R$ and $q\in Q$. Since $ab\in Q$, $aK\subseteq Q$; whence $ay\in Q$. Now, $a^{k+1}z=a(a^{k}z)=a(y-q)\in Q$. This implies that $z\in A_{k+1}$, and so  $sz\in sA_{k+1}\subseteq A_k$. Consequently, $a^{k}sz\in Q $ which implies that  $a^{k}sz +sq=sy\in Q$. Thus we have $s(I\cap K)\subseteq Q\subseteq I\cap K$. This implies that there exists $s'\in S$ such that either $ss'I\subseteq Q$ or $ss'K\subseteq Q$ since $Q$ is  $S$-irreducible. If $ss'K\subseteq Q$, then $ss'b\in Q$  which is not possible. Therefore $ss'I\subseteq Q$ which implies that $ss'a^{k}\in Q$. Put $t=ss'\in S$. Then $(ta)^{k} \in Q$ and hence  $ta\in rad(Q)$, as desired.
\end{proof}

\begin{theorem}\label{on}
\noindent
Let $R$ be an $S$-$J$-Noetherian ring. Then every proper $J$-ideal of $R$ disjoint with $S$ can be written as a finite intersection of $S$-primary ideals.
\end{theorem}

\begin{proof}
	\noindent
	Let $E$ be the collection of $J$-ideals of $R$ which are disjoint with $S$ and can not be written as a finite intersection of $S$-primary ideals. We wish to show $E=\emptyset$. On the contrary suppose  $E\neq\emptyset$.  Since $R$ is an $S$-$J$-Noetherian ring, by Theorem \ref{j1}, there exists an $S$-maximal element in $E$, say $I$. Evidently, $I$ is not an $S$-primary ideal, by Lemma \ref{sp}, $I$ is not an $S$-irreducible ideal, and so  $I$ is not an irreducible ideal. This implies that  $I=K\cap L$ for some ideals $K$ and $L$ of $R$ with $I\neq K$ and $I\neq L$. As $I$ is not $S$-irreducible, and so $sK\nsubseteq I$ and $sL\nsubseteq I$ for all $s\in S$. Now, we claim that $K$, $L\notin E$. For this, if $K$ (respectively, $L$) belongs to $E$, then since $I$ is an $S$-maximal element of $E$ and $I\subset K$ (respectively, $I\subset L$), there exists $s'$ (respectively, $s''$) from $S$ such that $s'K\subseteq I$ (respectively, $s''L\subseteq I$). This is not possible, as $I$ is not $S$-irreducible. Therefore $K, L\notin E$. Also, if $K\cap S\neq\emptyset$ $(respectively, L\cap S\neq\emptyset)$, then there exist $s_{1}\in K\cap S$  $(respectively, s_{2}\in L\cap S)$. This implies that $s's_{1}\in s'K\subseteq I$ $(respectively, s''s_{2}\in s''L\subseteq I)$, which is a contradiction because $I$ disjoint with $S$. Thus $K$ and $L$ are also disjoint with $S$. This implies that $K$ and $L$ can be written as a finite intersection of $S$-primary ideals. Consequently, $I$ can also be written as a finite intersection of $S$-primary ideals since  $I=K\cap L$, a contradiction as $I\in E$. Thus $E=\emptyset$, i.e., every proper $J$-ideal of $R$ disjoint with $S$ can be written as a finite intersection of $S$-primary ideals. 
\end{proof}

\end{document}